\documentclass[12pt,reqno]{amsart}
\usepackage{amsmath, amsfonts, amssymb, amsthm, hyperref}
\usepackage{bm}
\allowdisplaybreaks[4]
\def\l{\left}
\def\r{\right}
\def\bg{\bigg}
\def\({\bg(}
\def\){\bg)}
\def\t{\text}
\def\f{\frac}

\def\per{{\rm per}}

\def\ls{\le}
\def\gs{\geqslant}

\def\bi{\binom}
\def\al{\alpha}

\def\ga{\gamma}

\def\eq{\equiv}

\def\da{\delta}

\def\Proof{\noindent{\it Proof}}

\def\Z{\mathbb Z}
\def\C{\mathbb C}
\def\N{\mathbb N}

\def\1{{\bf 1}}

\def\jacob #1#2{\genfrac{(}{)}{}{}{#1}{#2}}
\def\pmod #1{\ ({\rm{mod}}\ #1)}

\def\<{\langle}
\def\>{\rangle}
\theoremstyle{plain}
\newtheorem{theorem}{Theorem}[section]
\newtheorem{lemma}{Lemma}
\newtheorem{corollary}{Corollary}

\theoremstyle{definition}

\theoremstyle{remark}
\newtheorem{remark}{Remark}

\makeatletter
\@namedef{subjclassname@2020}{%
  \textup{2020} Mathematics Subject Classification}
\makeatother
 \vspace{4mm}

\begin{document}
\hbox{Preprint}
\medskip

\title[Evaluations of some Toeplitz-type determinants]
{Evaluations of some Toeplitz-type determinants}
\author{Han Wang}
\address {(Han Wang) Department of Mathematics, Nanjing
University, Nanjing 210093, People's Republic of China}
\email{hWang@smail.nju.edu.cn}

\author{Zhi-Wei Sun}
\address {(Zhi-Wei Sun, corresponding author) Department of Mathematics, Nanjing
University, Nanjing 210093, People's Republic of China}
\email{zwsun@nju.edu.cn}

\keywords{Determinant, eigenvalue, Lucas sequence, Toeplitz matrix.
\newline \indent Supported by the National Natural Science Foundation of China (grant no. 11971222).}
\begin{abstract}
In this paper we evaluate some Toeplitz-type determinants.
Let $n>1$ be an integer. We prove the following two basic identities:
\begin{align*}
\det{[j-k+\delta_{jk}]_{1\leq j,k\leq n}}&=1+\frac{n^2(n^2-1)}{12},
\\
\det{[|j-k|+\delta_{jk}]_{1\leq j,k\leq n}}&=
\begin{cases}
\frac{1+(-1)^{(n-1)/2}n}{2}&\text{if}\ 2\nmid n,\\
\frac{1+(-1)^{n/2}}{2}&\text{if}\ 2\mid n,
\end{cases}
\end{align*}
where $\delta_{jk}$ is the Kronecker delta.
For complex numbers $a,b,c$ with $b\not=0$ and $a^2\not=4b$, and the sequence
$(w_m)_{m\in\mathbb Z}$ with $w_{k+1}=aw_k-bw_{k-1}$ for all $k\in\mathbb Z$, we establish the identity
$$\det[w_{j-k}+c\delta_{jk}]_{1\ls j,k\ls n}
=c^n+c^{n-1}nw_0+c^{n-2}(w_1^2-aw_0w_1+bw_0^2)\frac{u_n^2b^{1-n}-n^2}{a^2-4b},$$
where $u_0=0$, $u_1=1$ and $u_{k+1}=au_k-bu_{k-1}$ for all $k=1,2,\ldots$.
\end{abstract}
\maketitle

\section{Introduction}
\setcounter{lemma}{0}
\setcounter{theorem}{0}
\setcounter{equation}{0}
\setcounter{conjecture}{0}
\setcounter{remark}{0}
\setcounter{corollary}{0}

For a matrix $M=[a_{jk}]_{1\ls j,k\ls n}$ over the field $\C$ of complex numbers, we use $\det(M)$
or $\det[a_{jk}]_{1\ls j,k\ls n}$ to denote its determinant.
A Toeplitz matrix over $\C$ has the form $[a_{j-k}]_{1\ls j,k\ls n}$.
In this paper we evaluate determinants of several Toeplitz matrices.

In 1934 the evaluation of $\det[|j-k|]_{1\ls j,k\ls n}$ was proposed by R. Robinson as a problem in Amer. Math. Monthly, later its solutions appeared in \cite{RS}. As a result,
$$\det[|j-k|]_{1\ls j,k\ls n}=(-1)^{n-1}(n-1)2^{n-2}.$$
Moreover, the inverse of the matrix $[|j-k|]_{1\ls j,k\ls n}$
was found by M. Fiedler (cf. J. Todd \cite{Todd}).

In part (i) of our first theorem we determine the value of $\det[|j-k|+\da_{jk}]_{1\ls j,k\ls n}$,
 where the Kronecker $\da_{jk}$ is $1$ or $0$ according as $j=k$ or not.

\begin{theorem}\label{2} Let $n$ be any positive integer.

{\rm (i)} We have
\begin{equation}\label{|j-k|}
\det[|j-k|+\da_{jk}]_{1\ls j,k\ls n}=\begin{cases}
\f{1+(-1)^{(n-1)/2}n}{2}&\t{if}\ 2\nmid n,\\
\f{1+(-1)^{n/2}}{2}&\t{if}\ 2\mid n.
\end{cases}\end{equation}

{\rm (ii)} We have
\begin{equation}\label{F_|j-k|}
\det[F_{|j-k|}+\da_{jk}]_{1\le j,k\le n}=\begin{cases}
1&\t{if}\ n\equiv0,\pm1\pmod6,\\
0&\t{otherwise},
\end{cases}
\end{equation}
where the Fibonacci numbers $F_0,F_1,F_2,\ldots$ are defined by $F_0=0,\ F_1=1$ and the recurrence
$$F_{i+1}=F_i+F_{i-1}\ \ (i=1,2,3,\ldots).$$

{\rm (iii)} We have
\begin{equation}\label{|j-k,3|}
\det\l[\jacob{|j-k|}3-\da_{jk}\r]_{1\le j,k\le n}=\begin{cases}
1&\t{if}\ n\equiv0\pmod6,\\
-1&\t{if}\ n\equiv\pm1\pmod6,\\
0&\t{otherwise},
\end{cases}
\end{equation}
where $(\f{\cdot}3)$ denotes the Legendre symbol.
\end{theorem}

\begin{remark} The three parts of Theorem \ref{2} are somewhat similar.
For $n=1, 2, 3, \cdots$, let $f(n)$ denote the right-hand side of \eqref{|j-k|}.
The sequence
\[
(f(n))_{n\geq 1}=(1,0,-1,3,0,-3,1,5,0,-5,1,7,0,-7,1,\ldots)
\]
appeared as \cite[A166445]{OEIS} which was generated by P. Barry in 2009 as the Hankel transform $(\det[a_{j+k-1}]_{1\leq j,k\leq n})_{n\geq1}$ of the integer sequence
\[
(a_1, a_2, a_3, a_4, a_5, a_6, a_7, a_8, a_9, a_{10}, \ldots)=(1,0,0,1,2,4,8,17,38,88, \ldots)
\]
satisfying the recurrence
\[
a_n=\sum_{k=1}^{n-1}a_ka_{n-k}\ \ (n=5,6,7,\ldots)
\]
which has a combinatorial interpretation (cf. \cite[A025276]{OEIS}).
\end{remark}

Recall that the $q$-analogue of an integer $m$ is given by
$$[m]_q=\f{q^m-1}{q-1}.$$
Note that $\lim_{q\to1}[m]_q=m$.
A natural extension of Theorem 1.1(i) is to evaluate the determinant
$A_n(q)=[[|j-k|]_q+\da_{jk}]_{1\ls j,k\ls n}$. However, after serious attempts via {\tt Mathematica},
we could not find a general pattern for the exact value of $\det(A_n(q))$.

For any integer $n\gs3$, we have
$$\det[j-k]_{1\ls j,k\ls n}=0$$
since $(1-k)+(3-k)=2(2-k)$ for all $k=1,\ldots,n$.
However, it is nontrivial to evaluate the determinant $\det[j-k+\da_{jk}]_{1\ls j,k\ls n}$.
The sequence $\det[j-k+\da_{jk}]_{1\ls j,k\ls n}$ $(n=1,2,3,\ldots)$
appeared as \cite[A079034]{OEIS} which was generated by B. Cloitre in 2003,
its initial fifteen terms are
$$1,\, 2,\, 7,\, 21,\, 51,\, 106,\, 197,\, 337,\, 541,\, 826,\, 1211,\, 1717,\, 2367,\, 3186,\, 4201.$$
In 2013 C. Baker added a comment to \cite[A079034]{OEIS} which asserts that
$$\det[j-k+\da_{jk}]_{1\ls j,k\ls n}=1+\f{n^2(n^2-1)}{12}$$
without any proof or linked reference. It seems that Baker found this by guessing the recurrence
of the sequence via using the Maple package {\tt gfun}.

Let $n$ be any positive integer, and let $B_n$
be the $n\times n$ matrix $[j-k+\da_{jk}]_{1\leq j, k\leq n}$.
We observe that when $n>2$ the $n$ eigenvalues of $B_n$ are
$$\lambda_1=1+\f{n\sqrt{n^2-1}}{2\sqrt3}\,i,\ \ \lambda_2=1-\f{n\sqrt{n^2-1}}{2\sqrt3}\,i,\ \ \lambda_3=\cdots=\lambda_n=1.$$
This yields that
\begin{equation}\label{det-A}
\det(B_n)=1+\f{n^2(n^2-1)}{12}.
\end{equation}
The referee informs that this can be simply proved via the Weinstein-Aronszajn identity
(cf. \cite{WA}); instead of proving \eqref{det-A} in details, we will deduce a further extension of \eqref{det-A}.
Concerning the permanent of $B_n$,  motivated by \cite[Conj. 11.23]{S-book}
we conjecture that
$$\per(B_{p-1})\eq3\pmod p\ \ \t{and}\ \ \per(B_p)\eq 1+4p\pmod {p^2}.$$
for any odd prime $p$.

Theorems 1.2-1.4 in this paper are concerned with evaluations of some Toeplitz-type determinants involving exponential functions.
For a function $g:\Z\to\C$ with $g(|n|)$ growing exponentially
(where $\C$ denotes the field of complex numbers), we find that the Weinstein-Aronszajn identity is
impractical to evaluate $\det[g(j-k)+\da_{j,k}]_{1\ls j,k\ls n}$ since we could not find explicit
patterns for $n\times2$ matrix $X$ and $2\times n$ matrix $Y$ satisfying
$$\det[g(j-k)+\da_{jk}]_{1\ls j,k\ls n}=\det[I_n+XY],$$
where $I_n$ is the identity matrix of order $n$.

Our following result is a further generalization of \eqref{det-A}.

\begin{theorem}\label{Th5} Let $n\gs2$ be an integer, and let $q\not=0,1$ be a complex number. Then
the characteristic polynomial of the matrix $P=[q^{j-k}+t]_{1\le j,k\le n}$
is
\begin{equation}\label{qxt}\det(xI_n-P)=x^{n-2}(x^2-n(t+1)x+t(n^2-q^{1-n}[n]_q^2)).
\end{equation}
\end{theorem}

Putting $t=-1$ and replacing $x$ by $(q-1)x$ in Theorem \ref{Th5}, we immediately obtain the following corollary.

\begin{corollary} Let $n\gs2$ be an integer, and let $q\not=0,1$ be a complex number. For the
matrix $P_q=[[j-k]_q]_{1\le j,k\le n}$ we have
\begin{equation}\det[xI_n-P_q]=x^{n}+\f{q^{1-n}[n]_q^2-n^2}{(q-1)^2}x^{n-2}.
\end{equation}
\end{corollary}

Observe that
\begin{align*}\lim_{q\to1}\f{q^{1-n}[n]_q^2-n^2}{(q-1)^2}
&=\lim_{t\to0}\f{(t+1)^{1-n}(((t+1)^n-1)/t)^2-n^2}{t^2}
\\&=\lim_{t\to0}\f{(t+1)^{1-n}((\sum_{k=1}^n\bi nkt^{k-1})^2-n^2)+((t+1)^{1-n}-1)n^2}{t^2}
\\&=\lim_{t\to0}\(\f{(n+\bi n2t+\bi n3t^2+\cdots)^2-n^2}{(t+1)^{n-1}t^2}+n^2\f{1-(t+1)^{n-1}}{(t+1)^{n-1}t^2}\)
\\&=\bi n2^2+2n\bi n3+\lim_{t\to0}\l(2n\bi n2\f{t^{-1}}{(t+1)^{n-1}}-n^2\f{\sum_{k=1}^n\bi{n-1}kt^{k-2}}{(t+1)^{n-1}}\r)
\\&=\bi n2^2+2n\bi n3-n^2\bi{n-1}2
=\f{n^2(n^2-1)}{12}.
\end{align*}
So, Corollary 1.1 implies \eqref{det-A} for any integer $n\gs2$.

Applying Corollary 1.1 with $q=-1$, we find that
$$\det(xI_n-P_{-1})=x^n+\f{(-1)^{n-1}[n]_{-1}^2-n^2}4x^{n-2}$$
for any integer $n\gs2$. In particular,
\begin{equation}\det\l[\f{1-(-1)^{j-k}}2+\da_{j,k}\r]_{1\ls j,k\ls n}=\f{9-(-1)^n-2n^2}8.\end{equation}

Applying Theorem \ref{Th5} with $(t,x)=(-1,-2),(1,-1)$, we obtain the following result.

\begin{corollary} For any positive integer $n$, we have
\begin{equation}\det[2^{j-k}-1+2\da_{jk}]_{1\ls j,k\ls n}=\f{4^n-2^{n-1}n^2+1}2.
\end{equation}
and
\begin{equation}\det[2^{j-k}+1+\da_{j,k}]_{1\ls j,k\ls n}=(n+1)^2-2^{1-n}(2^n-1)^2.
\end{equation}
\end{corollary}

For complex numbers $a$ and $b\not=0$, the Lucas sequence $u_m=u_m(a,b)$ $(m\in\Z)$
and its companion sequence $v_m=v_m(a,b)$ $(m\in\Z)$
are defined as follows:
\begin{gather*}u_0=0,\ u_1=1,\ \t{and} \ u_{k+1}=au_k-bu_{k-1}\ \t{for all}\ k\in\Z;
\\v_0=2,\ v_1=a,\ \t{and} \ v_{k+1}=av_k-bv_{k-1}\ \t{for all}\ k\in\Z.
\end{gather*}
By the Binet formula,
$$(\al-\beta)u_m=\al^m-\beta^m\ \ \t{and}\ \ v_m=\al^m+\beta^m\ \ \t{for all}\ m\in\Z,$$
where
\begin{equation}\label{alpha}\al=\f{a+\sqrt{a^2-4b}}2\ \ \t{and}\ \ \beta=\f{a-\sqrt{a^2-4b}}2\end{equation}
are the two roots of the quadratic equation $x^2-ax+b=0$.
Clearly $b^nu_{-n}=-u_n$ and $b^nv_{-n}=v_n$ for all $n\in\N$.
 For any positive integer $n$, it is known that
$$u_n=\sum_{k=0}^{\lfloor(n-1)/2\rfloor}\bi{n-1-k}ka^{n-1-2k}(-b)^k
\ \ \ \t{and}\ \ \
v_n=\sum_{k=0}^{\lfloor n/2\rfloor}\f n{n-k}\bi{n-k}ka^{n-2k}(-b)^k,$$
which can be easily proved by induction.
Note also that $u_m(2,1)=m$ for all $m\in\Z$.

For $P(z)=\sum_{k=0}^{n-1}a_kz^k\in\C[z]$,  it is known (cf. \cite[Lemma 9]{K05}) that
$$\det[P(x_j+y_k)]_{1\ls j,k\ls n}=a_{n-1}^n\prod_{r=0}^{n-1}\bi{n-1}r\times\prod_{1\ls j<k\ls n}(x_j-x_k)(y_k-y_j).$$
Thus, for any integer $n\gs3$, and complex numbers $a$ and $b\not=0$, we have
\begin{equation}\label{uv}(\al-\beta)^n\det[u_{j-k}(a,b)]_{1\ls j,k\ls n}=\det[v_{j-k}(a,b)]_{1\ls j,k\ls n}=0
\end{equation}
(where $\al$ and $\beta$ are given by \eqref{alpha}),
since
\begin{align*}\det\l[\al^{j-k}\pm\beta^{j-k}\r]_{1\ls j,k\ls n}
=\prod_{k=1}^n\al^{-k}\times\prod_{j=1}^n\beta^j\times\det\l[\l(\f{\al}{\beta}\r)^j
\pm\l(\f{\al}{\beta}\r)^k\r]_{1\ls j,k\ls n}=0.
\end{align*}

In contrast with \eqref{uv}, we obtain the following new result.

\begin{theorem} \label{Th-w} Let $a$ and $b\not=0$ be complex numbers with $a^2\not=4b$.
Let $(w_m)_{m\in\Z}$ be a sequence of complex numbers with $w_{k+1}=aw_k-bw_{k-1}$ for all $k\in\Z$.
For any complex number $c$ and integer $n\gs2$, we have
\begin{equation}\label{w-value}\det[w_{j-k}+c\da_{jk}]_{1\ls j,k\ls n}
=c^n+c^{n-1}nw_0+c^{n-2}(w_1^2-aw_0w_1+bw_0^2)\f{u_n(a,b)^2b^{1-n}-n^2}{a^2-4b}.
\end{equation}
\end{theorem}

It is hard to guess the exact formula for $\det[w_{j-k}+c\da_{jk}]_{1\ls j,k\ls n}$
in Theorem \ref{Th-w}
via looking at various numerical examples.
The main trick for Theorem \ref{Th-w}
is how we find the identity \eqref{w-value} via eigenvalues and eigenvectors, this will be seen in Section 3.

\begin{corollary} \label{Th2} Let $a,b,c$ be complex numbers with $b\not=0$ and $a^2\not=4b$.
For any integer $n\gs2$, we have
\begin{equation}\label{A_2} \det[u_{j-k}(a,b)+c\da_{jk}]_{1\ls j,k\ls n}=c^n+c^{n-2}\f{u_n(a,b)^2b^{1-n}-n^2}{a^2-4b}.
\end{equation}
and
\begin{equation}\label{A_3} \det[v_{j-k}(a,b)+c\da_{jk}]_{1\ls j,k\ls n}=c^{n-2}((n+c)^2-u_n(a,b)^2b^{1-n}).
\end{equation}
\end{corollary}

For any $m\in\Z$, $u_m(-1,1)$ coincides with the Legendre symbol $(\f m3)$,
and $v_m(1,-1)=\omega^m+\bar\omega^m$ where $\omega$ denotes the cubic root $(-1+\sqrt{-3})/2$ of unity.
Applying Corollary \ref{Th2} with $a=-1$ and $b=1$, we get the following result.

\begin{corollary} For any integer $n\gs2$ and complex number $c$, we have
\begin{equation} \det\l[\l(\f{j-k}3\r)+c\da_{j,k}\r]_{1\ls j,k\ls n}=c^n+c^{n-2}\l\lfloor\f{n^2}3\r\rfloor,
\end{equation}
where $(\f{\cdot}3)$ is the Legendre symbol.
\end{corollary}

Recall that those $F_m=u_m(1,-1)$ ($m\in\Z$) are the well-known Fibonacci numbers, and those $L_m=v_m(1,-1)$ $(m\in\Z$)
are the Lucas numbers.
Corollary \ref{Th2} with $a=1$ and $b=-1$ yields the following result.

\begin{corollary} For any integer $n\gs2$ and complex number $c$, we have
\begin{equation}\label{F}\det[F_{j-k}+c\da_{jk}]_{1\ls j,k\ls n}=c^n+\f{c^{n-2}}5\l((-1)^{n-1}F_n^2-n^2\r)
\end{equation}
and
\begin{equation}\det[L_{j-k}+c\da_{jk}]_{1\ls j,k\ls n}=c^{n-2}((n+c)^2+(-1)^nF_n^2).
\end{equation}
\end{corollary}

In contrast with Theorem \ref{Th5}, we also establish the following result.

\begin{theorem}\label{Th6}
Let $n\gs2$ be an integer, and let $q\not=0,1$ be a complex number. For the matrix $Q=[q^{j+k}+t]_{0\le j,k\le n-1}$, we have
\begin{equation}\label{qxt2}\det(xI_n-Q)=x^n-(nt+[n]_{q^2})x^{n-1}+(n[n]_{q^2}-[n]_q^2)tx^{n-2}.
\end{equation}
\end{theorem}

The identity \eqref{qxt2} with $q=2$ and $x=t=-1$ yields the following corollary.

\begin{corollary} For any positive integer $n$, we have
\begin{equation}\det[2^{j+k}-1+\da_{jk}]_{0\ls j,k\ls n-1}=(2^n-1)^2-(n-1)\f{4^n+2}3.
\end{equation}
\end{corollary}

Although we have Theorem \ref{Th6} which is similar to Theorem \ref{Th5},
we don't have a result similar to Theorem \ref{Th-w}.

We are going to prove Theorem \ref{2} and Theorems \ref{Th5}-\ref{Th-w} in Sections 2 and 3, respectively.
Section 4 is devoted to our proof of Theorem \ref{Th6}.

\section{Proof of Theorem $\ref{2}$}
\setcounter{lemma}{0}
\setcounter{theorem}{0}
\setcounter{equation}{0}
\setcounter{conjecture}{0}
\setcounter{remark}{0}
\setcounter{corollary}{0}

Let $n$ be a positive integer. For $j_0,k_0\in\{1,\ldots,n\}$ with $j_0\not=k_0$, we define
$$T_{j_0, k_0}=[t_{jk}]_{1\leq j,k\leq n}\ \ \t{and}\ \ \widetilde T_{j_0, k_0}=[\widetilde t_{jk}]_{1\leq j,k\leq n},$$
where
\[
t_{jk}=\begin{cases}
1&\t{if}\ j=k, \\
-1&\t{if}\ j=j_0\ \t{and}\ k=k_0, \\
0&\t{otherwise},
\end{cases}
\]
and
\[
\widetilde t_{jk}=\begin{cases}
1&\t{if}\ j=j_0\ \t{and}\ k=k_0, \t{or}\ j=k,\\
0&\t{otherwise}.
\end{cases}
\]
It is easy to see that
\begin{equation}\label{TT}\det(T_{j_0,k_0})=1=\det(\widetilde T_{j_0,k_0})\quad\t{for all}\ j_0,k_0=1,\ldots,n\ \t{with}\ j_0\not=k_0.
\end{equation}
 We need this useful fact in our proofs of Theorem \ref{2}.

\medskip
\noindent {\it Proof of Theorem} \ref{2}(i). Let $A$ denote the matrix $[|j-k|+\da_{jk}]_{1\ls j,k\ls n}$.
Clearly \eqref{|j-k|} holds trivially for $n=1,2$.
(When $n=2$ all the entries of $A$ are $1$.)

Now we assume that $n\gs3$.
Observe that
\[
T_{21}T_{32}\cdots T_{n-1,n-2}T_{n,n-1}AT_{n-1,n}T_{n-2,n-1}\cdots T_{23}T_{12}=C,
\]
where $C=[c_{jk}]_{1\leq j,k\leq n}$ and
\[
c_{jk}=\begin{cases}
1&\t{if}\ 1\in\{j,k\}\ \t{and}\ jk\neq2, \\
-1&\t{if}\ |j-k|=1\ \t{and}\ jk\neq 2, \\
0&\t{otherwise}.
\end{cases}
\]
It follows that $\det(A)=\det(C)$ in view of \eqref{TT}.

Note that both the last column and the last row of $T_{n,n-2}CT_{n-2,n}$ contain
a unique nonzero entry (which is $1$). We illustrate this via the transformation from $C$ to $T_{n,n-2}CT_{n-2,n}$:
\begin{equation}\label{m}
\begin{bmatrix}
1&0&1&1&\cdots&1&1&1&1\\
0&0&-1&0&\cdots&0&0&0&0\\
1&-1&0&-1&\cdots&0&0&0&0\\
1&0&-1&0&\cdots&0&0&0&0\\
\vdots&\vdots&\vdots&\vdots&\ddots&\vdots&\vdots&\vdots\\
1&0&0&0&\cdots&0&-1&0&0\\
1&0&0&0&\cdots&-1&0&-1&0\\
1&0&0&0&\cdots&0&-1&0&-1\\
1&0&0&0&\cdots&0&0&-1&0
\end{bmatrix}
\to
\begin{bmatrix}
1&0&1&1&\cdots&1&1&1&0\\
0&0&-1&0&\cdots&0&0&0&0\\
1&-1&0&-1&\cdots&0&0&0&0\\
1&0&-1&0&\cdots&0&0&0&0\\
\vdots&\vdots&\vdots&\vdots&\ddots&\vdots&\vdots&\vdots\\
1&0&0&0&\cdots&0&-1&0&1\\
1&0&0&0&\cdots&-1&0&-1&0\\
1&0&0&0&\cdots&0&-1&0&0\\
0&0&0&0&\cdots&1&0&0&0
\end{bmatrix}.
\end{equation}
Thus, via expanding $\det(T_{n,n-2}CT_{n-2,n})$ by its last column and the last row, we obtain
$$\det(A)=\det(C)=\det(T_{n,n-2}CT_{n-2,n})=-\det{D_{n-2}},$$
where the $(n-2)\times(n-2)$ matrix $D_{n-2}$ has the form
\[
\begin{bmatrix}
1&0&1&1&\cdots&1&1&1&1\\
0&0&-1&0&\cdots&0&0&0&0\\
1&-1&0&-1&\cdots&0&0&0&0\\
1&0&-1&0&\cdots&0&0&0&0\\
\vdots&\vdots&\vdots&\vdots&\ddots&\vdots&\vdots&\vdots\\
1&0&0&0&\cdots&0&-1&0&0\\
1&0&0&0&\cdots&-1&0&0&0\\
1&0&0&0&\cdots&0&0&0&-1\\
1&0&0&0&\cdots&0&0&-1&0
\end{bmatrix}
\]
obtained by deleting the $n$th and the $(n-3)$th columns and rows from the last matrix in $\eqref{m}$.

By repeating the above procedure with $T_{n,n-2}CT_{n-2,n}\to D_{n-2}$ replaced by
\[
T_{n-4k,n-2-4k}D_{n-2k}T_{n-2-4k,n-4k}\to D_{n-2-2k}\ \ \l(0<k<\l\lfloor\f n4\r\rfloor\r),
\]
we get that
\begin{align*}\det(A)&=-\det(D_{n-2})=\det(D_{n-4})=\cdots
\\&= (-1)^{\lfloor n/4\rfloor}\det\l(D_{n-2\lfloor\f n4\rfloor}\r)=(-1)^{\lfloor n/4\rfloor}\det\l(D_{\f{n+\{n\}_4}2}\r)
\end{align*}
with the aid of \eqref{TT},
where $\{n\}_4$ is the least nonnegative residue of $n$ modulo $4$.
\medskip

{\it Case} 1. $n\eq0\pmod2$.

When $n\equiv0\pmod4$, the matrix $D_{\f{n+\{n\}_4}{2}}=D_{\f n2}$ has the form
\[
\begin{bmatrix}
0&-1&0&0&\cdots&0&0&0&0\\
-1&0&0&0&\cdots&0&0&0&0\\
0&0&0&-1&\cdots&0&0&0&0\\
0&0&-1&0&\cdots&0&0&0&0\\
\vdots&\vdots&\vdots&\vdots&\ddots&\vdots&\vdots&\vdots\\
0&0&0&0&\cdots&0&-1&0&0\\
0&0&0&0&\cdots&-1&0&0&0\\
0&0&0&0&\cdots&0&0&0&-1\\
0&0&0&0&\cdots&0&0&-1&0
\end{bmatrix},
\]
therefore
\[
\det\l(D_\f{n}{2}\r)=(-1)^{n/2}\times(-1)^{n/4}=(-1)^{n/4}
\]
and hence $
\det(A)=(-1)^{n/4}\det(D_\f{n}{2})=1.$

When $n\equiv2\pmod4$, the matrix $D_{\f{n+\{n\}_4}2}=D_{\f n2+1}$
has the form
\[
\begin{bmatrix}
1&0&1&1&\cdots&1&1&1&1\\
0&0&0&0&\cdots&0&0&0&0\\
1&0&0&-1&\cdots&0&0&0&0\\
1&0&-1&0&\cdots&0&0&0&0\\
\vdots&\vdots&\vdots&\vdots&\ddots&\vdots&\vdots&\vdots\\
1&0&0&0&\cdots&0&-1&0&0\\
1&0&0&0&\cdots&-1&0&0&0\\
1&0&0&0&\cdots&0&0&0&-1\\
1&0&0&0&\cdots&0&0&-1&0
\end{bmatrix}
\]
(with the second row containing only zero entries),
and hence
\[
\det(A)=(-1)^{(n-2)/4}\det\l(D_{\f{n}{2}+1}\r)=0.
\]

Below we assume that $n$ is odd.
\medskip

{\it Case} 2. $n\equiv1\pmod4$.

In this case, we have
\[D_{\f{n+\{n\}_4}2}=D_\f{n+1}{2}
=\begin{bmatrix}
1&1&1&1&1&\cdots&1&1\\
1&0&-1&0&0&\cdots&0&0\\
1&-1&0&0&0&\cdots&0&0\\
1&0&0&0&-1&\cdots&0&0\\
1&0&0&-1&0&\cdots&0&0\\
\vdots&\vdots&\vdots&\vdots&\vdots&\ddots&\vdots&\vdots\\
1&0&0&0&0&\cdots&0&-1\\
1&0&0&0&0&\cdots&-1&0
\end{bmatrix},
\]
 and
 $M_\f{n+1}{2}=\widetilde T_{1,\f{n+1}{2}}\widetilde T_{1,\f{n-1}{2}}\cdots \widetilde T_{12}D_\f{n+1}{2}$ has the form
\[
\begin{bmatrix}
\f{n+1}{2}&0&0&0&0&\cdots&0&0\\
1&0&-1&0&0&\cdots&0&0\\
1&-1&0&0&0&\cdots&0&0\\
1&0&0&0&-1&\cdots&0&0\\
1&0&0&-1&0&\cdots&0&0\\
\vdots&\vdots&\vdots&\vdots&\vdots&\ddots&\vdots&\vdots\\
1&0&0&0&0&\cdots&0&-1\\
1&0&0&0&0&\cdots&-1&0
\end{bmatrix}.
\]
Thus
\[
\det\l(D_\f{n+1}{2}\r)=\det\l(M_\f{n+1}{2}\r)=\f{n+1}2(-1)^{(n-1)/2}\times(-1)^{(n-1)/4}=\f{n+1}2(-1)^{(n-1)/4}
\]
and hence
\[
\det(A)=(-1)^{(n-1)/4}\det\l(D_\f{n+1}{2}\r)=\f{n+1}{2}.
\]

{\it Case} 3. $n\equiv3\pmod4$.

In this case, we have
\[D_{\f{n+\{n\}_4}2}=D_\f{n+3}{2}=
\begin{bmatrix}
1&0&1&1&1&\cdots&1&1\\
0&0&-1&0&0&\cdots&0&0\\
1&-1&0&0&0&\cdots&0&0\\
1&0&0&0&-1&\cdots&0&0\\
1&0&0&-1&0&\cdots&0&0\\
\vdots&\vdots&\vdots&\vdots&\vdots&\ddots&\vdots&\vdots\\
1&0&0&0&0&\cdots&0&-1\\
1&0&0&0&0&\cdots&-1&0
\end{bmatrix},
\]
and $M_\f{n+3}{2}=\widetilde T_{1,\f{n+3}{2}}\cdots\widetilde T_{12}D_\f{n+1}{2}\widetilde T_{21}\cdots\widetilde  T_{\f{n+3}{2},1}$ has the form
\[
\begin{bmatrix}
\f{n-1}{2}&-1&0&0&0&\cdots&0&0\\
-1&0&-1&0&0&\cdots&0&0\\
0&-1&0&0&0&\cdots&0&0\\
0&0&0&0&-1&\cdots&0&0\\
0&0&0&-1&0&\cdots&0&0\\
\vdots&\vdots&\vdots&\vdots&\vdots&\ddots&\vdots&\vdots\\
0&0&0&0&0&\cdots&0&-1\\
0&0&0&0&0&\cdots&-1&0
\end{bmatrix}.
\]
Therefore
\[
\det\l(D_\f{n+3}{2}\r)=\det\l(M_\f{n+3}{2}\r)=\f{n-1}2(-1)^{(n+1)/2}\times(-1)^{(n+1)/4}=\f{n-1}2(-1)^{(n+1)/4}
\]
and hence
\[
\det(A)=(-1)^{(n-3)/4}\det\l(D_\f{n+3}{2}\r)=\f{1-n}{2}.
\]

In view of the above, we have completed our proof of the first part of Theorem \ref{2}.

\medskip
\noindent {\it Proof of Theorem} \ref{2}(ii). Let $A_n$ denote the matrix $[F_{|j-k|}+\da_{jk}]_{1\le j,k\le n}$.
It is easy to verify \eqref{F_|j-k|} for $1\le n\le 6$.

Now, let $n\ge7$. It suffices to prove that $\det(A_n)=\det(A_{n-6})$.
Note that  \[A_n=\begin{bmatrix}
1&1&1&2&\cdots&F_{n-4}&F_{n-3}&F_{n-2}&F_{n-1}\\
1&1&1&1&\cdots&F_{n-5}&F_{n-4}&F_{n-3}&F_{n-2}\\
1&1&1&1&\cdots&F_{n-6}&F_{n-5}&F_{n-4}&F_{n-3}\\
2&1&1&1&\cdots&F_{n-7}&F_{n-6}&F_{n-5}&F_{n-4}\\
\vdots&\vdots&\vdots&\vdots&\ddots&\vdots&\vdots&\vdots&\vdots\\
F_{n-4}&F_{n-5}&F_{n-6}&F_{n-7}&\cdots&1&1&1&2\\
F_{n-3}&F_{n-4}&F_{n-5}&F_{n-6}&\cdots&1&1&1&1\\
F_{n-2}&F_{n-3}&F_{n-4}&F_{n-5}&\cdots&1&1&1&1\\
F_{n-1}&F_{n-2}&F_{n-3}&F_{n-4}&\cdots&2&1&1&1
\end{bmatrix}. \]
Since
\[T_{n,n-2}T_{n,n-1}A_nT_{n-1,n}T_{n-2,n}=\begin{bmatrix}
1&1&1&2&\cdots&F_{n-4}&F_{n-3}&F_{n-2}&0\\
1&1&1&1&\cdots&F_{n-5}&F_{n-4}&F_{n-3}&0\\
1&1&1&1&\cdots&F_{n-6}&F_{n-5}&F_{n-4}&0\\
2&1&1&1&\cdots&F_{n-7}&F_{n-6}&F_{n-5}&0\\
\vdots&\vdots&\vdots&\vdots&\ddots&\vdots&\vdots&\vdots&\vdots\\
F_{n-4}&F_{n-5}&F_{n-6}&F_{n-7}&\cdots&1&1&1&0\\
F_{n-3}&F_{n-4}&F_{n-5}&F_{n-6}&\cdots&1&1&1&-1\\
F_{n-2}&F_{n-3}&F_{n-4}&F_{n-5}&\cdots&1&1&1&-1\\
0&0&0&0&\cdots&0&-1&-1&1
\end{bmatrix}, \]
 the matrix $\widetilde T_{n-2,n}\widetilde T_{n-1,n}T_{n,n-2}T_{n,n-1}A_nT_{n-1,n}T_{n-2,n}\widetilde T_{n,n-1}\widetilde T_{n,n-2}$ has the form
\[\begin{bmatrix}
1&1&1&2&\cdots&F_{n-4}&F_{n-3}&F_{n-2}&0\\
1&1&1&1&\cdots&F_{n-5}&F_{n-4}&F_{n-3}&0\\
1&1&1&1&\cdots&F_{n-6}&F_{n-5}&F_{n-4}&0\\
2&1&1&1&\cdots&F_{n-7}&F_{n-6}&F_{n-5}&0\\
\vdots&\vdots&\vdots&\vdots&\ddots&\vdots&\vdots&\vdots&\vdots\\
F_{n-4}&F_{n-5}&F_{n-6}&F_{n-7}&\cdots&1&1&1&0\\
F_{n-3}&F_{n-4}&F_{n-5}&F_{n-6}&\cdots&1&0&0&0\\
F_{n-2}&F_{n-3}&F_{n-4}&F_{n-5}&\cdots&1&0&0&0\\
0&0&0&0&\cdots&0&0&0&1
\end{bmatrix}. \]
As there is a unique nonzero entry (which is $1$) in the last row of the last matrix, and
$$\det\l(\widetilde T_{n-2,n}\widetilde T_{n-1,n}T_{n,n-2}T_{n,n-1}A_nT_{n-1,n}T_{n-2,n}\widetilde T_{n,n-1}\widetilde T_{n,n-2}\r)=\det(A_n)$$
in light of \eqref{TT}, we see that $\det(A_n)=\det(A^{(1)})$, where
\[A^{(1)}=\begin{bmatrix}
1&1&1&2&\cdots&F_{n-5}&F_{n-4}&F_{n-3}&F_{n-2}\\
1&1&1&1&\cdots&F_{n-6}&F_{n-5}&F_{n-4}&F_{n-3}\\
1&1&1&1&\cdots&F_{n-7}&F_{n-6}&F_{n-5}&F_{n-4}\\
2&1&1&1&\cdots&F_{n-8}&F_{n-7}&F_{n-6}&F_{n-5}\\
\vdots&\vdots&\vdots&\vdots&\ddots&\vdots&\vdots&\vdots&\vdots\\
F_{n-5}&F_{n-6}&F_{n-7}&F_{n-8}&\cdots&1&1&1&2\\
F_{n-4}&F_{n-5}&F_{n-6}&F_{n-7}&\cdots&1&1&1&1\\
F_{n-3}&F_{n-4}&F_{n-5}&F_{n-6}&\cdots&1&1&0&0\\
F_{n-2}&F_{n-3}&F_{n-4}&F_{n-5}&\cdots&2&1&0&0
\end{bmatrix}.
\]
Similarly, $\widetilde T_{n-3,n-1}\widetilde T_{n-2,n-1}T_{n-1,n-3}T_{n-1,n-2}A^{(1)}T_{n-2,n-1}T_{n-3,n-1}\widetilde T_{n-1,n-2}\widetilde T_{n-1,n-3}$ has the form
\[\begin{bmatrix}
1&1&1&2&\cdots&F_{n-5}&F_{n-4}&F_{n-3}&0\\
1&1&1&1&\cdots&F_{n-6}&F_{n-5}&F_{n-4}&0\\
1&1&1&1&\cdots&F_{n-7}&F_{n-6}&F_{n-5}&0\\
2&1&1&1&\cdots&F_{n-8}&F_{n-7}&F_{n-6}&0\\
\vdots&\vdots&\vdots&\vdots&\ddots&\vdots&\vdots&\vdots&\vdots\\
F_{n-5}&F_{n-6}&F_{n-7}&F_{n-8}&\cdots&1&1&1&0\\
F_{n-4}&F_{n-5}&F_{n-6}&F_{n-7}&\cdots&1&0&0&0\\
F_{n-3}&F_{n-4}&F_{n-5}&F_{n-6}&\cdots&1&0&-1&0\\
0&0&0&0&\cdots&0&0&0&1
\end{bmatrix},
\]
and hence $\det(A^{(1)})$ equals the determinant of the matrix
\[A^{(2)}=\begin{bmatrix}
1&1&1&2&\cdots&F_{n-6}&F_{n-5}&F_{n-4}&F_{n-3}\\
1&1&1&1&\cdots&F_{n-7}&F_{n-6}&F_{n-5}&F_{n-4}\\
1&1&1&1&\cdots&F_{n-8}&F_{n-7}&F_{n-6}&F_{n-5}\\
2&1&1&1&\cdots&F_{n-9}&F_{n-8}&F_{n-7}&F_{n-6}\\
\vdots&\vdots&\vdots&\vdots&\ddots&\vdots&\vdots&\vdots&\vdots\\
F_{n-6}&F_{n-7}&F_{n-8}&F_{n-9}&\cdots&1&1&1&2\\
F_{n-5}&F_{n-6}&F_{n-7}&F_{n-8}&\cdots&1&1&1&1\\
F_{n-4}&F_{n-5}&F_{n-6}&F_{n-7}&\cdots&1&1&0&0\\
F_{n-3}&F_{n-4}&F_{n-5}&F_{n-6}&\cdots&2&1&0&-1
\end{bmatrix}.
\]
Note also that
\begin{align*}&\widetilde T_{n-4,n-2}\widetilde T_{n-3,n-2}T_{n-2,n-4}T_{n-2,n-3}A^{(2)}T_{n-3,n-2}T_{n-4,n-2}
\\=\ &\begin{bmatrix}
1&1&1&2&\cdots&F_{n-6}&F_{n-5}&F_{n-4}&0\\
1&1&1&1&\cdots&F_{n-7}&F_{n-6}&F_{n-5}&0\\
1&1&1&1&\cdots&F_{n-8}&F_{n-7}&F_{n-6}&0\\
2&1&1&1&\cdots&F_{n-9}&F_{n-8}&F_{n-7}&0\\
\vdots&\vdots&\vdots&\vdots&\ddots&\vdots&\vdots&\vdots&\vdots\\
F_{n-6}&F_{n-7}&F_{n-8}&F_{n-9}&\cdots&1&1&1&0\\
F_{n-5}&F_{n-6}&F_{n-7}&F_{n-8}&\cdots&1&0&0&-1\\
F_{n-4}&F_{n-5}&F_{n-6}&F_{n-7}&\cdots&1&0&-1&-1\\
0&0&0&0&\cdots&0&-1&-1&0
\end{bmatrix}
\end{align*}
and hence
$$T_{n-3,n-5}T_{n-3,n-4}\widetilde T_{n-4,n-2}\widetilde T_{n-3,n-2}T_{n-2,n-4}T_{n-2,n-3}A^{(2)}T_{n-3,n-2}T_{n-4,n-2}T_{n-4,n-3}T_{n-5,n-3}$$ has the form
\[\begin{bmatrix}
1&1&1&2&\cdots&F_{n-6}&F_{n-5}&0&0\\
1&1&1&1&\cdots&F_{n-7}&F_{n-6}&0&0\\
1&1&1&1&\cdots&F_{n-8}&F_{n-7}&0&0\\
2&1&1&1&\cdots&F_{n-9}&F_{n-8}&0&0\\
\vdots&\vdots&\vdots&\vdots&\ddots&\vdots&\vdots&\vdots&\vdots\\
F_{n-6}&F_{n-7}&F_{n-8}&F_{n-9}&\cdots&1&1&-1&0\\
F_{n-5}&F_{n-6}&F_{n-7}&F_{n-8}&\cdots&1&0&-1&-1\\
0&0&0&0&\cdots&-1&-1&0&0\\
0&0&0&0&\cdots&0&-1&0&0
\end{bmatrix}.
\]
Note that both the last row and the last column contain only one nonzero entry (which is $-1$). Therefore, $\det(A^{(2)})$ equals the determinant of the matrix
\[\begin{bmatrix}
1&1&1&2&\cdots&F_{n-7}&F_{n-6}&0&0&0\\
1&1&1&1&\cdots&F_{n-8}&F_{n-7}&0&0&0\\
1&1&1&1&\cdots&F_{n-9}&F_{n-8}&0&0&0\\
2&1&1&1&\cdots&F_{n-10}&F_{n-9}&0&0&0\\
\vdots&\vdots&\vdots&\vdots&\ddots&\vdots&\vdots&\vdots&\vdots&\vdots\\
F_{n-7}&F_{n-8}&F_{n-9}&F_{n-10}&\cdots&1&1&0&0&0\\
F_{n-6}&F_{n-7}&F_{n-8}&F_{n-9}&\cdots&1&1&0&-1&0\\
0&0&0&0&\cdots&0&0&0&0&-1\\
0&0&0&0&\cdots&0&-1&0&0&0\\
0&0&0&0&\cdots&0&0&-1&0&0
\end{bmatrix},
\]
which has a unique nonzero entry (which is $-1$) in the $(n-3)$-th row and the $(n-3)$-th column.
Thus, $\det(A^{(2)})$ coincides with the determinant of the matrix
\[A_{n-6}=\begin{bmatrix}
1&1&1&2&\cdots&F_{n-10}&F_{n-9}&F_{n-8}&F_{n-7}\\
1&1&1&1&\cdots&F_{n-11}&F_{n-10}&F_{n-9}&F_{n-8}\\
1&1&1&1&\cdots&F_{n-12}&F_{n-11}&F_{n-10}&F_{n-9}\\
2&1&1&1&\cdots&F_{n-13}&F_{n-12}&F_{n-11}&F_{n-10}\\
\vdots&\vdots&\vdots&\vdots&\ddots&\vdots&\vdots&\vdots&\vdots\\
F_{n-10}&F_{n-11}&F_{n-12}&F_{n-13}&\cdots&1&1&1&2\\
F_{n-9}&F_{n-10}&F_{n-11}&F_{n-12}&\cdots&1&1&1&1\\
F_{n-8}&F_{n-9}&F_{n-10}&F_{n-11}&\cdots&1&1&1&1\\
F_{n-7}&F_{n-8}&F_{n-9}&F_{n-10}&\cdots&2&1&1&1
\end{bmatrix}.
\]

In view of the above, $\det(A_n)=\det(A^{(1)})=\det(A^{(2)})=\det(A_{n-6})$.
This concludes our proof of Theorem \ref2(ii). \qed

\medskip
\noindent{\it Proof of Theorem} \ref2(iii). Let $A_n$ denote the matrix $[(\f{|j-k|}3)-\da_{jk}]_{1\le j,k\le n}$.
It is easy to verify \eqref{|j-k,3|} for $1\le n\le 6$.

Now, let $n\ge7$. It suffices to prove that $\det(A_n)=\det(A_{n-6})$.
Note that \[A_n=\begin{bmatrix}
-1&1&-1&0&\cdots&\jacob{n-4}3&\jacob{n-3}3&\jacob{n-2}3&\jacob{n-1}3\\
1&-1&1&-1&\cdots&\jacob{n-5}3&\jacob{n-4}3&\jacob{n-3}3&\jacob{n-2}3\\
-1&1&-1&1&\cdots&\jacob{n-6}3&\jacob{n-5}3&\jacob{n-4}3&\jacob{n-3}3\\
0&-1&1&-1&\cdots&\jacob{n-7}3&\jacob{n-6}3&\jacob{n-5}3&\jacob{n-4}3\\
\vdots&\vdots&\vdots&\vdots&\ddots&\vdots&\vdots&\vdots&\vdots\\
\jacob{n-4}3&\jacob{n-5}3&\jacob{n-6}3&\jacob{n-7}3&\cdots&-1&1&-1&0\\
\jacob{n-3}3&\jacob{n-4}3&\jacob{n-5}3&\jacob{n-6}3&\cdots&1&-1&1&-1\\
\jacob{n-2}3&\jacob{n-3}3&\jacob{n-4}3&\jacob{n-5}3&\cdots&-1&1&-1&1\\
\jacob{n-1}3&\jacob{n-2}3&\jacob{n-3}3&\jacob{n-4}3&\cdots&0&-1&1&-1
\end{bmatrix}. \]
Also,
\begin{align*} &\widetilde T_{n,n-2}\widetilde T_{n,n-1}A_n\widetilde T_{n-1,n}\widetilde T_{n-2,n}
=\begin{bmatrix}
-1&1&-1&0&\cdots&\jacob{n-4}3&\jacob{n-3}3&\jacob{n-2}3&0\\
1&-1&1&-1&\cdots&\jacob{n-5}3&\jacob{n-4}3&\jacob{n-3}3&0\\
-1&1&-1&1&\cdots&\jacob{n-6}3&\jacob{n-5}3&\jacob{n-4}3&0\\
0&-1&1&-1&\cdots&\jacob{n-7}3&\jacob{n-6}3&\jacob{n-5}3&0\\
\vdots&\vdots&\vdots&\vdots&\ddots&\vdots&\vdots&\vdots&\vdots\\
\jacob{n-4}3&\jacob{n-5}3&\jacob{n-6}3&\jacob{n-7}3&\cdots&-1&1&-1&0\\
\jacob{n-3}3&\jacob{n-4}3&\jacob{n-5}3&\jacob{n-6}3&\cdots&1&-1&1&-1\\
\jacob{n-2}3&\jacob{n-3}3&\jacob{n-4}3&\jacob{n-5}3&\cdots&-1&1&-1&1\\
0&0&0&0&\cdots&0&-1&1&-1
\end{bmatrix}, \end{align*}
and hence $T_{n-2,n}\widetilde T_{n-1,n}\widetilde T_{n,n-2}\widetilde T_{n,n-1}A_n\widetilde T_{n-1,n}\widetilde T_{n-2,n}\widetilde T_{n,n-1}T_{n,n-2}$ has the form
\[\begin{bmatrix}
-1&1&-1&0&\cdots&\jacob{n-4}3&\jacob{n-3}3&\jacob{n-2}3&0\\
1&-1&1&-1&\cdots&\jacob{n-5}3&\jacob{n-4}3&\jacob{n-3}3&0\\
-1&1&-1&1&\cdots&\jacob{n-6}3&\jacob{n-5}3&\jacob{n-4}3&0\\
0&-1&1&-1&\cdots&\jacob{n-7}3&\jacob{n-6}3&\jacob{n-5}3&0\\
\vdots&\vdots&\vdots&\vdots&\ddots&\vdots&\vdots&\vdots&\vdots\\
\jacob{n-4}3&\jacob{n-5}3&\jacob{n-6}3&\jacob{n-7}3&\cdots&-1&1&-1&0\\
\jacob{n-3}3&\jacob{n-4}3&\jacob{n-5}3&\jacob{n-6}3&\cdots&1&0&0&0\\
\jacob{n-2}3&\jacob{n-3}3&\jacob{n-4}3&\jacob{n-5}3&\cdots&-1&0&0&0\\
0&0&0&0&\cdots&0&0&0&-1
\end{bmatrix}. \]
As there is a unique nonzero entry (which is $-1$) in the last row,
and $$\det\l(T_{n-2,n}\widetilde T_{n-1,n}\widetilde T_{n,n-2}\widetilde T_{n,n-1}A_n\widetilde T_{n-1,n}\widetilde T_{n-2,n}\widetilde T_{n,n-1}T_{n,n-2}\r)=\det(A_n)$$
in light of \eqref{TT},
we have $-\det(A_n)=\det(A^{(1)})$, where
\[A^{(1)}=\begin{bmatrix}
-1&1&-1&0&\cdots&\jacob{n-5}3&\jacob{n-4}3&\jacob{n-3}3&\jacob{n-2}3\\
1&-1&1&-1&\cdots&\jacob{n-6}3&\jacob{n-5}3&\jacob{n-4}3&\jacob{n-3}3\\
-1&1&-1&1&\cdots&\jacob{n-7}3&\jacob{n-6}3&\jacob{n-5}3&\jacob{n-4}3\\
0&-1&1&-1&\cdots&\jacob{n-8}3&\jacob{n-7}3&\jacob{n-6}3&\jacob{n-5}3\\
\vdots&\vdots&\vdots&\vdots&\ddots&\vdots&\vdots&\vdots&\vdots\\
\jacob{n-5}3&\jacob{n-6}3&\jacob{n-7}3&\jacob{n-8}3&\cdots&-1&1&-1&0\\
\jacob{n-4}3&\jacob{n-5}3&\jacob{n-6}3&\jacob{n-7}3&\cdots&1&-1&1&-1\\
\jacob{n-3}3&\jacob{n-4}3&\jacob{n-5}3&\jacob{n-6}3&\cdots&-1&1&0&0\\
\jacob{n-2}3&\jacob{n-3}3&\jacob{n-4}3&\jacob{n-5}3&\cdots&0&-1&0&0
\end{bmatrix}.
\]
Similarly, $T_{n-3,n-1}\widetilde T_{n-2,n-1}\widetilde T_{n-1,n-3}\widetilde T_{n-1,n-2}A^{(1)}\widetilde T_{n-2,n-1}\widetilde T_{n-3,n-1}\widetilde T_{n-1,n-2}T_{n-1,n-3}$ has the form
\[\begin{bmatrix}
-1&1&-1&0&\cdots&\jacob{n-5}3&\jacob{n-4}3&\jacob{n-3}3&0\\
1&-1&1&-1&\cdots&\jacob{n-6}3&\jacob{n-5}3&\jacob{n-4}3&0\\
-1&1&-1&1&\cdots&\jacob{n-7}3&\jacob{n-6}3&\jacob{n-5}3&0\\
0&-1&1&-1&\cdots&\jacob{n-8}3&\jacob{n-7}3&\jacob{n-6}3&0\\
\vdots&\vdots&\vdots&\vdots&\ddots&\vdots&\vdots&\vdots&\vdots\\
\jacob{n-5}3&\jacob{n-6}3&\jacob{n-7}3&\jacob{n-8}3&\cdots&-1&1&-1&0\\
\jacob{n-4}3&\jacob{n-5}3&\jacob{n-6}3&\jacob{n-7}3&\cdots&1&0&0&0\\
\jacob{n-3}3&\jacob{n-4}3&\jacob{n-5}3&\jacob{n-6}3&\cdots&-1&0&1&0\\
0&0&0&0&\cdots&0&0&0&-1
\end{bmatrix},
\]
and hence $\det(A^{(1)})=-\det(A^{(2)})$, where
\[A^{(2)}=\begin{bmatrix}
-1&1&-1&0&\cdots&\jacob{n-6}3&\jacob{n-5}3&\jacob{n-4}3&\jacob{n-3}3\\
1&-1&1&-1&\cdots&\jacob{n-7}3&\jacob{n-6}3&\jacob{n-5}3&\jacob{n-4}3\\
-1&1&-1&1&\cdots&\jacob{n-8}3&\jacob{n-7}3&\jacob{n-6}3&\jacob{n-5}3\\
0&-1&1&-1&\cdots&\jacob{n-9}3&\jacob{n-8}3&\jacob{n-7}3&\jacob{n-6}3\\
\vdots&\vdots&\vdots&\vdots&\ddots&\vdots&\vdots&\vdots&\vdots\\
\jacob{n-6}3&\jacob{n-7}3&\jacob{n-8}3&\jacob{n-9}3&\cdots&-1&1&-1&0\\
\jacob{n-5}3&\jacob{n-6}3&\jacob{n-7}3&\jacob{n-8}3&\cdots&1&-1&1&-1\\
\jacob{n-4}3&\jacob{n-5}3&\jacob{n-6}3&\jacob{n-7}3&\cdots&-1&1&0&0\\
\jacob{n-3}3&\jacob{n-4}3&\jacob{n-5}3&\jacob{n-6}3&\cdots&0&-1&0&1
\end{bmatrix}.
\]
Also,
\begin{align*}&T_{n-4,n-2}\widetilde T_{n-3,n-2}\widetilde T_{n-2,n-4}\widetilde T_{n-2,n-3}A^{(2)}\widetilde T_{n-3,n-2}\widetilde T_{n-4,n-2}
\\=\ &\begin{bmatrix}
-1&1&-1&0&\cdots&\jacob{n-6}3&\jacob{n-5}3&\jacob{n-4}3&0\\
1&-1&1&-1&\cdots&\jacob{n-7}3&\jacob{n-6}3&\jacob{n-5}3&0\\
-1&1&-1&1&\cdots&\jacob{n-8}3&\jacob{n-7}3&\jacob{n-6}3&0\\
0&-1&1&-1&\cdots&\jacob{n-9}3&\jacob{n-8}3&\jacob{n-7}3&0\\
\vdots&\vdots&\vdots&\vdots&\ddots&\vdots&\vdots&\vdots&\vdots\\
\jacob{n-6}3&\jacob{n-7}3&\jacob{n-8}3&\jacob{n-9}3&\cdots&-1&1&-1&0\\
\jacob{n-5}3&\jacob{n-6}3&\jacob{n-7}3&\jacob{n-8}3&\cdots&1&0&0&-1\\
\jacob{n-4}3&\jacob{n-5}3&\jacob{n-6}3&\jacob{n-7}3&\cdots&-1&0&1&1\\
0&0&0&0&\cdots&0&-1&1&0
\end{bmatrix},
\end{align*}
and hence $$T_{n-3,n-5}T_{n-3,n-4}T_{n-4,n-2}\widetilde T_{n-3,n-2}\widetilde T_{n-2,n-4}\widetilde T_{n-2,n-3}A^{(2)}\widetilde T_{n-3,n-2}\widetilde T_{n-4,n-2}T_{n-4,n-3}T_{n-5,n-3}$$ has the form
\[\begin{bmatrix}
-1&1&-1&0&\cdots&\jacob{n-6}3&\jacob{n-5}3&0&0\\
1&-1&1&-1&\cdots&\jacob{n-7}3&\jacob{n-6}3&0&0\\
-1&1&-1&1&\cdots&\jacob{n-8}3&\jacob{n-7}3&0&0\\
0&-1&1&-1&\cdots&\jacob{n-9}3&\jacob{n-8}3&0&0\\
\vdots&\vdots&\vdots&\vdots&\ddots&\vdots&\vdots&\vdots&\vdots\\
\jacob{n-6}3&\jacob{n-7}3&\jacob{n-8}3&\jacob{n-9}3&\cdots&-1&1&-1&0\\
\jacob{n-5}3&\jacob{n-6}3&\jacob{n-7}3&\jacob{n-8}3&\cdots&1&0&1&-1\\
0&0&0&0&\cdots&-1&1&0&0\\
0&0&0&0&\cdots&0&-1&0&0
\end{bmatrix}.
\]
 Since there is only one nonzero entry (which is $-1$) in the last row of the last matrix, we see that $\det(A^{(2)})$ equals the determinant of the matrix
\[\begin{bmatrix}
-1&1&-1&0&\cdots&\jacob{n-7}3&\jacob{n-6}3&0&0&0\\
1&-1&1&-1&\cdots&\jacob{n-8}3&\jacob{n-7}3&0&0&0\\
-1&1&-1&1&\cdots&\jacob{n-9}3&\jacob{n-8}3&0&0&0\\
0&-1&1&-1&\cdots&\jacob{n-10}3&\jacob{n-9}3&0&0&0\\
\vdots&\vdots&\vdots&\vdots&\ddots&\vdots&\vdots&\vdots&\vdots&\vdots\\
\jacob{n-7}3&\jacob{n-8}3&\jacob{n-9}3&\jacob{n-10}3&\cdots&-1&1&0&0&0\\
\jacob{n-6}3&\jacob{n-7}3&\jacob{n-8}3&\jacob{n-9}3&\cdots&1&-1&0&-1&0\\
0&0&0&0&\cdots&0&0&0&0&-1\\
0&0&0&0&\cdots&0&-1&0&0&0\\
0&0&0&0&\cdots&0&0&-1&0&0
\end{bmatrix},
\]
which has a unique nonzero entry (which is $-1$) in the $(n-3)$-th row and in the $(n-3)$-th column.
Thus, with the aid of \eqref{TT},  $\det(A^{(2)})$ equals the determinant of the matrix
\[A_{n-6}=\begin{bmatrix}
-1&1&-1&0&\cdots&\jacob{n-10}3&\jacob{n-9}3&\jacob{n-8}3&\jacob{n-7}3\\
1&-1&1&-1&\cdots&\jacob{n-11}3&\jacob{n-10}3&\jacob{n-9}3&\jacob{n-8}3\\
-1&1&-1&1&\cdots&\jacob{n-12}3&\jacob{n-11}3&\jacob{n-10}3&\jacob{n-9}3\\
0&-1&1&-1&\cdots&\jacob{n-13}3&\jacob{n-12}3&\jacob{n-11}3&\jacob{n-10}3\\
\vdots&\vdots&\vdots&\vdots&\ddots&\vdots&\vdots&\vdots&\vdots\\
\jacob{n-10}3&\jacob{n-11}3&\jacob{n-12}3&\jacob{n-13}3&\cdots&-1&1&-1&0\\
\jacob{n-9}3&\jacob{n-10}3&\jacob{n-11}3&\jacob{n-12}3&\cdots&1&-1&1&-1\\
\jacob{n-8}3&\jacob{n-9}3&\jacob{n-10}3&\jacob{n-11}3&\cdots&-1&1&-1&1\\
\jacob{n-7}3&\jacob{n-8}3&\jacob{n-9}3&\jacob{n-10}3&\cdots&0&-1&1&-1
\end{bmatrix}.
\]
Therefore, $\det(A_n)=-\det(A^{(1)})=\det(A^{(2)})=\det(A_{n-6})$. This ends our proof of Theorem \ref2(iii).
 \qed

\section{Proofs of Theorems \ref{Th5}-\ref{Th-w}}
\setcounter{lemma}{0}
\setcounter{theorem}{0}
\setcounter{equation}{0}
\setcounter{conjecture}{0}
\setcounter{remark}{0}
\setcounter{corollary}{0}

\begin{lemma} \label{Lem-x3} Let $n$ be a positive integer, and let $q\neq0$ and $t$ be complex numbers
with $n-[n]_q+t(q^{1-n}[n]_q-n)\not=0$. Suppose that
\begin{equation}\label{al-ga3} \ga=\f{n(t+1)\pm\sqrt{n^2(t-1)^2+4tq^{1-n}[n]_q^2}}2
\end{equation}
and
\begin{equation}\label{al-x3}
y=\f{\ga-[n]_q-nt}{n-[n]_q+(q^{1-n}[n]_q-n)t}.
\end{equation}
Then
\begin{equation}\label{da3}
\sum_{k=1}^n(q^{j-k}+t)(1+y(q^{k-n}-1))=\ga(1+y(q^{j-n}-1))
\end{equation}
for any positive integer $j$.
\end{lemma}
\Proof. As
$$\gamma^2-n(t+1)\gamma+(n^2-q^{1-n}[n]_q^2)t=0,$$
we have
$$[n]_q(n-[n]_q+(q^{1-n}[n]_q-n)t)=(\ga-[n]_q-nt)(\ga-n+[n]_q)$$
and hence
\begin{equation}\label{x3} (\ga-n+[n]_q)y=[n]_q.
\end{equation}

Let $j\in\{1,2,3,\ldots\}$, and set
$$\Delta_j=\sum_{k=1}^n(q^{j-k}+t)(1+y(q^{k-n}-1))-\ga(1+y(q^{j-n}-1)).$$
Then
\[\label{da}\begin{aligned}&\Delta_j-t(1+y(q^{k-n}-1))+\ga(1-y)
\\&=q^{j-n}\(\sum_{k=1}^nq^{n-k}(1+y(q^{k-n}-1))-\ga y\)
\\&=q^{j-n}\l([n]_q(1-y)+ny-\ga y)\r)=0
\end{aligned}\]
by \eqref{x3}. So $\Delta_1=\Delta_2=\cdots$.

Next we show that $\Delta_n=0$.
Observe that
\begin{align*}
&\sum_{k=1}^n(q^{n-k}+t)(1+(q^{k-n}-1)y)
\\=\ &\sum_{k=1}^n\l(q^{n-k}(1-y)+t(1-y)+y+q^{k-n}ty\r)
\\=\ &[n]_q(1-y)+nt(1-y)+ny+q^{1-n}[n]_qty
\\=\ &[n]_q+nt+y(n-[n]_q+(q^{1-n}[n]_q-n)t)
\\=\ &\ga=\ga(1+y(q^{n-n}-1))
\end{align*}
by the definition of $y$. So $\Delta_n=0$.

In view of the above, $\Delta_j=0$ for all $j=1,2,3,\ldots$. This concludes our proof. \qed

\medskip
\noindent{\it Proof of Theorem \ref{Th5}}. It is easy to verify the desired result for $n=2$.
Below we assume that $n\gs3$.

If $n-[n]_q$ and $q^{1-n}[n]_q-n$ are both zero, then $q^{n-1}=1$
and $n=[n]_q=1$. As $n\gs3$,  there are infinitely many $t\in\C$ such that
$$\begin{cases}n-[n]_q+t(q^{1-n}[n]_q-n)\not=0,
\\n^2(t-1)^2+4tq^{1-n}[n]_q^2\not=0.\end{cases}$$
Take any such a number $t$, and choose $\gamma$ and $y$
as in  \eqref{al-ga3} and \eqref{al-x3}. Then $\ga$ given by \eqref{al-ga3} is an eigenvalue of the matrix $P=[q^{j-k}+t]_{1\le j,k\le n}$, and the column vector
$v=(v_1,\ldots,v_n)^T$ with $v_k=1+y(q^{k-n}-1)$ is an eigenvector of $P$ associated with the eigenvalue $\ga$.
Note that $\gamma$ given by \eqref{al-ga3} has two different choices since $n^2(t-1)^2+4tq^{1-n}[n]_q^2\not=0$.

Let $s\in\{3,\ldots,n\}$. For $1\ls k\ls n$, let us define
$$v^{(s)}_k=
\begin{cases}
q^{2-s}[s-2]_q&\t{if}\ k=1, \\
-q^{2-s}[s-1]_q&\t{if}\ k=2, \\
\da_{sk}&\t{if}\ 3\ls k\ls n.
\end{cases}$$
It is easy to verify that
$$\sum_{k=1}^nv^{(s)}_k=0=\sum_{k=1}^nq^{j-k}v^{(s)}_k\ \ \ \t{for all}\ j=1,\ldots,n.$$
Thus $0$ is an eigenvalue of the matrix $P=[q^{j-k}+t]_{1\le j,k\le n}$, and
the column vector $v^{(s)}=(v^{(s)}_1,\ldots,v^{(s)}_n)^T$ is an eigenvector
of $P$ associated with the eigenvalue $0$.

If $\sum_{s=3}^n c_sv^{(s)}$ is the zero column vector for some $c_3,\ldots,c_n\in\C$,
then for each $k=3,\ldots,n$ we have
$$c_k=\sum_{s=3}^n c_s\da_{sk}=\sum_{s=3}^n c_sv^{(s)}_k=0.$$
Thus the $n-2$ column vectors $v^{(3)},\ldots,v^{(n)}$ are linearly independent over $\C$.

By the above, the $n$ eigenvalues of the matrix $P=[q^{j-k}+t]_{1\le j,k\le n}$ are the two values of $\ga$ given by \eqref{da3}, and $\lambda_3=\cdots=\lambda_n=0$.
Thus the characteristic polynomial of $P$
is
\begin{align*}\det(xI_n-P)=&\l(x-\f{n(t+1)}2-\f{\sqrt{n^2(t-1)^2+4tq^{1-n}[n]_q^2}}2\r)
\\&\times\l(x-\f{n(t+1)}2+\f{\sqrt{n^2(t-1)^2+4tq^{1-n}[n]_q^2}}2\r)\prod_{s=3}^n(x-\lambda_s)
\\=&x^{n-2}\l(\l(x-\f{n(t+1)}2\r)^2-\f{n^2(t-1)^2+4tq^{1-n}[n]_q^2}4\r)
\\=&x^{n-2}(x^2-n(t+1)x+t(n^2-q^{1-n}[n]_q^2)).
\end{align*}

In light of the above, the identity \eqref{qxt} holds for infinitely many values of $t$.
Note that both sides of \eqref{qxt} are polynomials in $t$ for any fixed $x\in\C$.
Thus, if we view both sides of \eqref{qxt} as polynomials in $x$ and $t$, then the identity
 \eqref{qxt} still holds. This ends our proof. \qed
 \medskip

\noindent{\it Proof of Theorem \ref{Th-w}}. If $w_0=w_1=0$ or $n=2$, then the desired result can be easily verified. Below we assume that $n\gs3$ and $\{w_0,w_1\}\not=\{0\}$.

(i) Let $\al$ and $\beta$ be the two roots of the quadratic equation
$z^2-az+b=0$. Note that $\al\beta=b\not=0$. Also, $\al\not=\beta$ since $\Delta=a^2-4b$ is nonzero.

It is well known that there are constants $c_1,c_2\in\C$ such that
$$w_m=c_1\al^m+c_2\beta^m\quad\t{for all}\ m\in\Z.$$
As $c_1+c_2=w_0$ and $c_1\al+c_2\beta=w_1$, we find that
\begin{equation}\label{cc} c_1=\f{w_1-\beta w_0}{\al-\beta}\ \ \t{and}\ \ c_2=\f{\al w_0-w_1}{\al-\beta}.\end{equation}
Since $w_0$ or $w_1$ is nonzero, one of $c_1$ and $c_2$ is nonzero.
Without any loss of generality, we assume $c_1\not=0$.

 Let $W$ denote the matrix $[w_{j-k}+c\da_{jk}]_{1\ls j,k\ls n}$. Then
\begin{align*}\det(W)=\ &\det\l[c_1\al^{j-k}+c_2\beta^{j-k}+c\da_{jk}\r]_{1\ls j,k\ls n}
\\=\ &c_1^n\prod_{j=1}^n\beta^{j}\times\prod_{k=1}^n\beta^{-k}
\times\det\l[\l(\f{\al}{\beta}\r)^{j-k}+\f{c_2+c\da_{jk}\beta^{k-j}}{c_1}\r]_{1\ls j,k\ls n}
\\=\ &c_1^n\det\l[q^{j-k}+t-x\da_{jk}\r]_{1\ls j,k\ls n}=(-c_1)^n\det\l[x\da_{jk}-q^{j-k}-t\r]_{1\ls j,k\ls n}
\end{align*}
where $q=\al/\beta\not=0,1$, and $t=c_2/c_1$ and $x=-c/c_1$.
Applying Theorem \ref{Th5}, we deduce that
\begin{align*}\det(W)&=(-c_1)^n x^{n-2}(x^2-n(t+1)x+t(n^2-q^{1-n}[n]_q^2))
\\&=c^{n-2}\l(c^2+nc(c_1+c_2)+c_1c_2\l(n^2-\f{\al^{1-n}}{\beta^{1-n}}\l(\f{(\al/\beta)^n-1}{\al/\beta-1}\r)^2\r)\r) \\&=c^n+nw_0c^{n-1}+c^{n-2}c_1c_2\l(n^2-(\al\beta)^{1-n}\l(\f{\al^n-\beta^n}{\al-\beta}\r)^2\r)
\\&=c^n+nw_0c^{n-1}+c^{n-2}c_1c_2\l(n^2-b^{1-n}u_n(a,b)^2\r).
\end{align*}
In view of \eqref{cc},
$$c_1c_2=\f{(w_1-\beta w_0)(\al w_0-w_1)}{(\al-\beta)^2}=\f{-w_1^2+(\al+\beta)w_0w_1-\al\beta w_0^2}{\Delta}=-\f{w_1^2-aw_0w_1+bw_0^2}{a^2-4b}.$$
Therefore the desired \eqref{w-value} follows. \qed

\section{Proof of Theorem \ref{Th6}}
\setcounter{lemma}{0}
\setcounter{theorem}{0}
\setcounter{equation}{0}
\setcounter{conjecture}{0}
\setcounter{remark}{0}
\setcounter{corollary}{0}

The following lemma is quite similar to Lemma \ref{Lem-x3}.

\begin{lemma} \label{Lem-x4} Let $n$ be a positive integer, and let $q\neq0$ and $t$ be complex numbers
with  $[n]_{q^2}+(q^{1-n}t-q^{n-1})[n]_q-nt\neq 0$. Suppose that
\begin{equation}\label{al-ga4}
\ga=\f{nt+[n]_{q^2}\pm\sqrt{(nt-[n]_{q^2})^2+4t[n]_q^2}}2
\end{equation}
and
\begin{equation}\label{al-x4}
z=\f{\ga-q^{n-1}[n]_q-nt}{[n]_{q^2}+(q^{1-n}t-q^{n-1})[n]_q-nt}.
\end{equation}
Then
\begin{equation}\label{da4}
\sum_{k=0}^{n-1}(q^{j+k}+t)(1+z(q^{k-n+1}-1))=\ga(1+z(q^{j-n+1}-1))
\end{equation}
for all $j=0,1,2,\ldots$.
\end{lemma}
\Proof. Since
$\ga^2-(nt+[n]_{q^2})\ga+t(n[n]_{q^2}-[n]_q^2)=0,$
we have
\begin{equation}\label{x4} (\ga-[n]_{q^2}+q^{n-1}[n]_q)z=q^{n-1}[n]_q.
\end{equation}

Let $j\in\{0,1,2,\ldots\}$, and set
$$R_j=\sum_{k=0}^{n-1}(q^{j+k}+t)(1+z(q^{k-n+1}-1))-\ga(1+z(q^{j-n+1}-1)).$$
It is easy to see that
\[\label{da}\begin{aligned} R_j-\sum_{k=0}^{n-1}t(1+z(q^{k-n+1}-1))+\ga(1-z)
=q^{j-n+1}\l({q^{n-1}[n]_q(1-z)+z[n]_{q^2}}-\ga z\r)=0
\end{aligned}\]
with the aid of \eqref{x4}. So $R_0=R_1=\cdots$. As
\begin{align*}
\sum_{k=0}^{n-1}(q^{n-1+k}+t)(1+z(q^{k-n+1}-1))
=\ga=\ga(1+z(q^{(n-1)-n+1}-1))
\end{align*}
we get $R_{n-1}=0$. So the desired result follows. \qed
\medskip

\noindent{\it Proof of Theorem \ref{Th6}}. It is easy to verify the desired result for $n=2$.
Below we assume that $n\gs3$.

If $[n]_{q^2}-q^{n-1}[n]_q$ and $q^{1-n}[n]_q-n$ are both zero, then
$[n]_q\not=0$ and
$$(q^n+1)[n]_q=(q+1)[n]_{q^2}=(q+1)q^{n-1}[n]_q=(q^n+q^{n-1})[n]_q,$$
hence $q^{n-1}=1$ and $n=[n]_q=1$.
 As $n\gs3$,  there are infinitely many $t\in\C$ such that
$$\begin{cases}[n]_{q^2}+(q^{1-n}t-q^{n-1})[n]_q-nt\not=0,
\\(nt-[n]_{q^2})^2+4t[n]_q^2\not=0.\end{cases}$$
Take any such a number $t$, and choose $\gamma$ and $z$
as in  \eqref{al-ga4} and \eqref{al-x4}. Then $\ga$ given by \eqref{al-ga4} is an eigenvalue of the matrix $Q=[q^{j+k}+t]_{0\le j,k\le n-1}$, and the column vector
$v=(v_0,\ldots,v_{n-1})^T$ with $v_k=1+z(q^{k-n+1}-1)$ is an eigenvector of $Q$ associated with the eigenvalue $\ga$.
Note that $\gamma$ given by \eqref{al-ga4} has two different choices since $(nt-[n]_{q^2})^2+4t[n]_q^2\not=0$.

Let $s\in\{3,\ldots,n\}$. For $k\in\{0,\ldots,n-1\}$, let us define
$$v^{(s)}_k=
\begin{cases}
q[s-2]_q&\t{if}\ k=0, \\
-[s-1]_q&\t{if}\ k=1, \\
\da_{s,k+1}&\t{if}\ 2\ls k\ls n-1.
\end{cases}$$
It is easy to verify that
$$\sum_{k=0}^{n-1}v^{(s)}_k=0=\sum_{k=0}^{n-1}q^{j+k}v^{(s)}_k\ \ \ \t{for all}\ j=1,\ldots,n.$$
Thus $0$ is an eigenvalue of the matrix $Q=[q^{j+k}+t]_{0\le j,k\le {n-1}}$, and
the column vector $v^{(s)}=(v^{(s)}_0,\ldots,v^{(s)}_{n-1})^T$ is an eigenvector
of $Q$ associated with the eigenvalue $0$.

If $\sum_{s=3}^{n} c_sv^{(s)}$ is the zero column vector for some $c_3,\ldots,c_{n}\in\C$,
then for each $k=2,\ldots,n-1$ we have
$$c_{k+1}=\sum_{s=3}^{n} c_s\da_{s,k+1}=\sum_{s=3}^{n} c_sv^{(s)}_{k}=0.$$
Thus the $n-2$ column vectors $v^{(3)},\ldots,v^{(n)}$ are linearly independent over $\C$.

By the above, the $n$ eigenvalues of the matrix $Q=[q^{j+k}+t]_{0\le j,k\le {n-1}}$ are the two values of $\ga$ given by \eqref{da4}, and $\lambda_3=\cdots=\lambda_n=0$.
Thus the characteristic polynomial of $Q$
is
\begin{align*}\det(xI_n-Q)=&\ \l(x-\f{nt+[n]_{q^2}}2-\f{\sqrt{(nt-[n]_{q^2})^2+4t[n]_q^2}}2\r)
\\&\ \times\l(x-\f{nt+[n]_{q^2}}2+\f{\sqrt{(nt-[n]_{q^2})^2+4t[n]_q^2}}2\r)\prod_{s=3}^n(x-\lambda_s)
\\=&\ x^{n-2}\l(\l(x-\f{nt+[n]_{q^2}}2\r)^2-\f{(nt-[n]_{q^2})^2+4t[n]_q^2}4\r)
\\=&\ x^n-(nt+[n]_{q^2})x^{n-1}+(n[n]_{q^2}-[n]_q^2)tx^{n-2}.
\end{align*}

In light of the above, the identity \eqref{qxt2} holds for infinitely many values of $t$.
Note that both sides of \eqref{qxt2} are polynomials in $t$ for any fixed $x\in\C$.
Thus, if we view both sides of \eqref{qxt2} as polynomials in $x$ and $t$, then the identity
 \eqref{qxt2} still holds. This concludes our proof. \qed

\medskip

\end{document}